\documentclass[11pt]{article}
\usepackage{amssymb}
\newtheorem{lemma}{Lemma}[section]
\newtheorem{definition}[lemma]{Definition}

\newtheorem{theorem}[lemma]{Theorem}
\newtheorem{proposition}[lemma]{Proposition}

\newtheorem{example}[lemma]{Example}

\newtheorem{note}[lemma]{Remark}

\def\dcl{\mathop{\rm dcl}\nolimits}
\def\Dom{\mathop{\rm Dom}\nolimits}
\def\Range{\mathop{\rm Range}\nolimits}

\title{Distributions of countable models \\ of quite
$o$-minimal Ehrenfeucht theories\footnote{This research was
partially supported by Committee of Science in Education and
Science Ministry of the Republic of Kazakhstan (Grant No.
AP05132546) and Russian Foundation for Basic Researches (Project No.
17-01-00531-a).}}
\author{B.Sh.~Kulpeshov, S.V.~Sudopalatov}
\date{}

\begin{document}

\maketitle

\begin{abstract}
We describe Rudin--Keisler preorders and distribution functions of
numbers of limit models for quite $o$-minimal Ehrenfeucht
theories. Decomposition formulas for these distributions are
found.

\end{abstract}

The notion of quite $o$-minimal was introduced and studied in
\cite{k2003}. This notion is a variation of weakly $o$-minimality
\cite{mms}. This notion occurred fruitful enough producing both
the structural description of models of these theories and the
generalization of Mayer theorem \cite{Ma}: it was shown that any
countable quite $o$-minimal theory has either finitely many
countable models, in the form $3^k\cdot 6^s$, or maximum many,
i.e. $2^\omega$, countable models \cite{KS}.

In the present paper, using a general theory of classification of
countable models of complete theories \cite{SuCCMCT1, SuCCMCT2} as
well as the description \cite{KS} of specificity for quite
$o$-minimal theories, we describe distributions of countable
models of quite $o$-minimal Ehrenfeucht theories in terms of
Rudin--Keisler preorders and distribution functions of numbers of
limit models. Besides, we derive decomposition formulas for these
distributions.

\section{Preliminaries}
\noindent

In this section we give the necessary information from
\cite{SuCCMCT1, SuCCMCT2}.

Recall that the number of pairwise non-isomorphic models of theory
$T$ and of cardinality $\lambda$ is denoted by
$I(T,\lambda)$\index{$I(T,\lambda)$}.

\begin{definition} {\rm \cite{Mi2} A theory $T$ is called
\emph{Ehrenfeucht} if $1<I(T,\omega)<\omega$.}
\end{definition}

\begin{definition} {\rm \cite{Be}  type $p(\bar{x})\in S(T)$ is said to be
\emph{powerful}\index{Type!powerful} in a theory $T$ if every
model ${\cal M}$ of $T$ realizing $p$ also realizes every type
$q\in S(T)$, that is, ${\cal M}\models S(T)$.}
\end{definition}

Since for any type $p\in S(T)$ there exists a countable model
${\cal M}$ of $T$, realizing $p$, and the model ${\cal M}$
realizes exactly countably many types, the availability of a
powerful type implies that $T$ is
\emph{small}\index{Theory!small}, that is, the set $S(T)$ is
countable. Hence for any type $q\in S(T)$ and its realization
$\bar{a}$, there exists a {\em prime model} ${\cal
M}({\bar{a}})$\index{${\cal M}({\bar{a}})$} {\em over}
$\bar{a}$,\index{Model!prime over a tuple} i.~e., a model of $T$
containing $\bar{a}$ with ${\cal M}(\bar{a})\models q(\bar{a})$
and such that ${\cal M}({\bar{a}})$ is elementarily embeddable to
any model realizing the type $q$. Since all prime models over
realizations of $q$ are isomorphic, we denote these models by
${\cal M}_q$\index{${\cal M}_q$}. Models ${\cal M}_q$ are called
{\em almost prime}\index{Model!almost prime} or {\em
$q$-prime}.\index{Model!$q$-prime}

\begin{definition} {\rm \cite{SuCCMCT1, Lr3, Su041} Let $p$ and $q$ be types in $S(T)$. We say that the type $p$
\emph{is dominated by a type}\index{Type!dominated} $q$, or $p$
\emph{does not exceed $q$ under the Rudin--Keisler
preorder}\index{Type!not exceed $q$ under the Rudin--Keisler
preorder}\index{Preorder!Rudin--Keisler} (written $p\leq_{\rm RK}
q$\index{$p\leq_{\rm RK}q$}\index{$\leq_{\rm RK}$}), if ${\cal
M}_q\models p$, that is, ${\cal M}_p$ is an elementary submodel of
${\cal M}_q$ (written ${\cal M}_p\preceq{\cal M}_q$). Besides, we
say that a~model ${\cal M}_p$ \emph{is dominated by a
model}\index{Model!dominated} ${\cal M}_q$, or ${\cal M}_p$
\emph{does not exceed ${\cal M}_q$ under the Rudin--Keisler
preorder}\index{Model!not exceed $q$ under the Rudin--Keisler
preorder}, and write ${\cal M}_p\leq_{\rm RK}{\cal
M}_q$\index{${\cal M}_p\leq_{\rm RK}{\cal M}_q$}.}
\end{definition}

Syntactically,  the  condition \ $p\leq_{\rm RK}q$  (and  hence
also ${\cal M}_p\leq_{\rm RK}{\cal M}_q$)  is expressed thus:
there exists a formula $\varphi(\bar{x},\bar{y})$ such that the
set $q(\bar{y})\cup\{\varphi(\bar{x},\bar{y})\}$ is consistent and
$q(\bar{y})\cup\{\varphi(\bar{x},\bar{y})\}\vdash p(\bar{x})$.
Since we deal with a small theory (there are only countably many
types over any tuple $\bar{a}$ and so any consistent formula with
parameters in $\bar{a}$ is deducible from a principal formula with
parameters in $\bar{a}$), $\varphi(\bar{x},\bar{y})$ can be chosen
so that for any formula $\psi(\bar{x},\bar{y})$, the set
$q(\bar{y})\cup\{\varphi(\bar{x},\bar{y}),\psi(\bar{x},\bar{y})\}$
being consistent implies that
$q(\bar{y})\cup\{\varphi(\bar{x},\bar{y})\}\vdash\psi(\bar{x},\bar{y})$.
In this event the formula $\varphi(\bar{x},\bar{y})$ is said to be
\emph{$(q,p)$-principal}\index{Formula!$(q,p)$-principal}.

\medskip
\begin{definition} {\rm \cite{SuCCMCT1, Lr3, Su041} Types $p$ and $q$ are said to be
\emph{domination-equivalent}\index{Types!domination-equivalent},
\emph{realization-equivalent}\index{Types!realization-equivalent},
\emph{Rudin--Keisler equivalent}\index{Types!Rudin--Keisler
equivalent}, or \ \emph{${\rm RK}$-equivalent}\index{Types!${\rm
RK}$-equivalent} \ (written \ $p\sim_{\rm RK} q$\index{$p\sim_{\rm
RK} q$}) if~$p\leq_{\rm RK} q$ and $q\leq_{\rm RK} p$. Models
${\cal M}_p$ and ${\cal M}_q$ are said to be
\emph{domination-equivalent}\index{Models!domination-equivalent},
\ \emph{Rudin--Keisler \ equivalent}\index{Models!Rudin--Keisler
equivalent}, \ or \ \emph{${\rm
RK}$-equivalent}\index{Models!${\rm RK}$-equivalent} \ (written \
\ ${\cal M}_p\sim_{\rm RK}{\cal M}_q$\index{${\cal M}_p\sim_{\rm
RK}{\cal M}_q$}).

As in \cite{Ta4}, types $p$ and $q$ are said to be \emph{strongly
domination-equivalent}\index{Types!domination-equivalent!strongly},
\emph{strongly
realization-equivalent}\index{Types!realization-equivalent!strongly},
\emph{strongly Rudin--Keisler
equivalent}\index{Types!Rudin--Keisler equivalent!strongly}, or
\emph{strongly ${\rm RK}$-equivalent}\index{Types!${\rm
RK}$-equivalent!strongly} (written \ $p\equiv_{\rm
RK}q$\index{$p\equiv_{\rm RK} q$}) if for some realizations
$\bar{a}$ and $\bar{b}$ of $p$ and $q$ respectively, both ${\rm
tp}(\bar{b}/\bar{a})$ and ${\rm tp}(\bar{a}/\bar{b})$ are
principal. Models ${\cal M}_p$ and ${\cal M}_q$ are said to be
\emph{strongly
domination-equivalent}\index{Models!domination-equivalent!strongly},
\emph{strongly Rudin--Keisler
equivalent}\index{Models!Rudin--Keisler equivalent}, or
\emph{strongly ${\rm RK}$-equivalent}\index{Models!${\rm
RK}$-equivalent!strongly} (written ${\cal M}_p\equiv_{\rm RK}{\cal
M}_q$\index{${\cal M}_p\equiv_{\rm RK}{\cal M}_q$}).}
\end{definition}

\medskip
Clearly, domination relations form preorders, and (strong)
do\-mi\-na\-tion-equi\-va\-lence \ relations \ are \ equivalence \
relations. \ Here, \ ${\cal M}_p\equiv_{\rm RK}{\cal M}_q$ implies
${\cal M}_p\sim_{\rm RK}{\cal M}_q$.

If ${\cal M}_p$ and ${\cal M}_q$ are not domina\-tion-equivalent
then they are non-isomorphic. Moreover, non-isomorphic models may
be found among domination-equivalent ones.

In \ Ehrenfeucht \ examples, \ models \ ${\cal
M}^n_{p_0},\ldots,{\cal M}^n_{p_{n-3}}$ \ are
domina\-tion-equivalent but pairwise non-isomorphic.

A syntactic characterization for the model isomorphism between
${\cal M}_p$ and ${\cal M}_q$ is given by the following
proposition. It asserts that the existence of an isomorphism
between ${\cal M}_p$ and ${\cal M}_q$ is equivalent to the strong
domination-equivalence of these models.

\medskip
\begin{proposition} {\rm \cite{SuCCMCT1, Su041}} For any types
$p(\bar{x})$ and $q(\bar{y})$ of a small theory $T$, the following
conditions are equivalent:

{\rm (1)} the models ${\cal M}_p$ and ${\cal M}_q$ are isomorphic;

{\rm (2)} the models ${\cal M}_p$ and ${\cal M}_q$ are strongly
domination-equivalent;

{\rm (3)} there exist $(p,q)$- and $(q,p)$-principal formulas
$\varphi_{p,q}(\bar{y},\bar{x})$ and
$\varphi_{q,p}(\bar{x},\bar{y})$ respectively, such that the set
$$
p(\bar{x})\cup q(\bar{y})\cup\{\varphi_{p,q}(\bar{y},\bar{x}),
\varphi_{q,p}(\bar{x},\bar{y})\}$$ is consistent;

{\rm (4)} there exists a $(p,q)$- and $(q,p)$-principal formula
$\varphi(\bar{x},\bar{y})$, such that the set
$$
p(\bar{x})\cup q(\bar{y})\cup\{\varphi(\bar{x},\bar{y})\}$$ is
consistent.
\end{proposition}

\begin{definition} {\rm \cite{SuCCMCT1, Su041}  Denote  by \ ${\rm RK}(T)$\index{${\rm RK}(T)$}  the  set
${\bf PM}$\index{${\bf PM}$}  of  isomorphism  types  of models
${\cal M}_p$, $p\in S(T)$, on which the relation of domination is
induced by $\leq_{\rm RK}$, a relation deciding domination among
${\cal M}_p$, that is, ${\rm RK}(T)=\langle{\bf PM};\leq_{\rm
RK}\rangle$. We say that isomorphism types ${\bf M}_1,{\bf
M}_2\in{\bf PM}$ are \emph{domination-equivalent} (written ${\bf
M}_1\sim_{\rm RK}{\bf M}_2$\index{${\bf M}_1\sim_{\rm RK}{\bf
M}_2$}) if so are their representatives.}
\end{definition}

Clearly, the preordered set ${\rm RK}(T)$ has a least element,
which is an isomorphism type of a prime model.

\medskip
\begin{proposition} {\rm \cite{SuCCMCT1, Su041}} If $I(T,\omega)<\omega$ then ${\rm
RK}(T)$ is a finite preordered set whose factor set ${\rm
RK}(T)/\!\!\!\sim_{\rm RK}$, with respect to
do\-mi\-na\-tion-equiva\-lence $\sim_{\rm RK}$, forms a partially
ordered set with a greatest element.
\end{proposition}

\begin{definition} {\rm \cite{SuCCMCT1, SuCCMCT2, Su041, Su08}
A model ${\cal M}$ of a theory $T$ is called {\em
limit}\index{Model!limit} if ${\cal M}$ is not prime over tuples
and ${\cal M}=\bigcup\limits_{n\in\omega}{\cal M}_n$ for some
elementary chain $({\cal M}_n)_{n\in\omega}$ of prime models of
$T$ over tuples. In this case the model ${\cal M}$ is said to be
{\em limit over a sequence ${\bf q}$ of
types}\index{Model!limit!over a sequence of types} or  {\em ${\bf
q}$-limit}\index{Model!${\bf q}$-limit}, where ${\bf
q}=(q_n)_{n\in\omega}$, ${\cal M}_n={\cal M}_{q_n}$, $n\in\omega$.
If the sequence ${\bf q}$ contains unique type $q$ then the ${\bf
q}$-limit model is called {\em limit over the type $q$}.}
\end{definition}

Denote \ by \ $I_p(T,\omega)$\index{$I_p(T,\omega)$} \ the \
number \ of \ pairwise \ non-isomorphic \ countable models of the
theory $T$, each of which is prime over a tuple, by $I_l(T)$ the
number of limit models of $T$, and by $I_l(T,q)$ the number of
limit models over a type $q\in S(T)$.

\begin{definition}
{\rm \cite{SuCCMCT2, Su08} A  theory  $T$  is  called  {\em
$p$-categorical}\index{Theory!$p$-categorical}  (respectively,
{\em $l$-categorical}\index{Theory!$l$-categorical},  {\em
$p$-Ehren\-feucht}\index{Theory!$p$-Ehrenfeucht}, and {\em
$l$-Ehrenfeucht})\index{Theory!$l$-Ehrenfeucht}  if
$I_p(T,\omega)=1$ \ (respectively, \ $I_l(T)=1$,
$1<I_p(T,\omega)<\omega$, $1<I_l(T)<\omega$).}
\end{definition}

Clearly, a small theory $T$ is $p$-categorical if and only if $T$
countably categorical, and if and only if $I_l(T)=0$; $T$ is
$p$-Ehrenfeucht if and only if the structure ${\rm RK}(T)$ finite
and has at least two elements; and $T$ is $p$-Ehrenfeucht with
$I_l(T)<\omega$ if and only if $T$ is Ehrenfeucht.

\medskip
Let $\widetilde{\bf M}\in{\rm RK}(T)/\!\!\sim_{\rm RK}$ be the
class consisting of isomorphism types of domination-equivalent
models ${\cal M}_{p_1},\ldots,{\cal M}_{p_n}$. Denote by ${\rm
IL}(\widetilde{\bf M})$ \index{${\rm IL}(\widetilde{\bf M})$} the
number of equivalence classes of models each of which is limit
over some type $p_i$.

\begin{theorem}\label{th1121_1136} {\rm \cite{SuCCMCT1, Su041}}
For any countable complete theory $T$, the following conditions
are equivalent:

{\rm (1)} $I(T,\omega)<\omega$;

{\rm (2)} $T$  is  small,  $|{\rm RK}(T)|<\omega$  and  ${\rm
IL}(\widetilde{\bf M})<\omega$  for  any  $\widetilde{\bf
M}\in{\rm RK}(T)/\!\!\sim_{\rm RK}$.

If $(1)$ or $(2)$ holds then $T$ possesses the following
properties:

{\rm (a)} ${\rm RK}(T)$ has a least element ${\bf M}_0$ {\rm (}an
isomorphism type of a prime model{\rm )} and ${\rm
IL}(\widetilde{{\bf M}_0})=0$;

{\rm (b)} ${\rm RK}(T)$ has a greatest $\sim_{\rm RK}$-class
$\widetilde{{\bf M}_1}$ {\rm (}a class of isomorphism types of all
prime models over realizations of powerful types{\rm )} and $|{\rm
RK}(T)|>1$ implies ${\rm IL}(\widetilde{{\bf M}_1})\geq 1$;

{\rm (c)} if $|\widetilde{\bf M}|>1$ then ${\rm IL}(\widetilde{\bf
M})\geq 1$.

Moreover, the following {\sl decomposition
formula}\index{Formula!decomposition} holds:
\begin{equation}\label{eqmain}
I(T,\omega)=|{\rm RK}(T)|+\sum_{i=0}^{|{\rm RK}(T)/\sim_{\rm
RK}|-1} {\rm IL}(\widetilde{{\bf M}_i}),
\end{equation}
where $\widetilde{{\bf M}_0},\ldots, \widetilde{{\bf M}_{|{\rm
RK}(T)/\sim_{\rm RK}|-1}}$ are all elements of the partially
ordered set ${\rm RK}(T)/\!\!\sim_{\rm RK}$.
\end{theorem}

\begin{definition} {\rm \cite{Wo} The {\em disjoint union}\index{Disjoint union!of
structures} $\bigsqcup\limits_{n\in\omega}{\cal
M}_n$\index{$\bigsqcup\limits_{n\in\omega}{\cal M}_n$} of pairwise
disjoint structures ${\cal M}_n$ for pairwise disjoint predicate
languages $\Sigma_n$, $n\in\omega$, is the structure of language
$\bigcup\limits_{n\in\omega}\Sigma_n\cup\{P^{(1)}_n\mid
n\in\omega\}$ with the universe $\bigsqcup\limits_{n\in\omega}
M_n$, $P_n=M_n$, and interpretations of predicate symbols in
$\Sigma_n$ coinciding with their interpretations in ${\cal M}_n$,
$n\in\omega$. The {\em disjoint union of theories}\index{Disjoint
union!of theories} $T_n$ for pairwise disjoint languages
$\Sigma_n$ accordingly, $n\in\omega$, is the theory
$$\bigsqcup\limits_{n\in\omega}T_n\rightleftharpoons{\rm Th}\left(\bigsqcup\limits_{n\in\omega}{\cal M}_n\right),$$
where\index{$\bigsqcup\limits_{n\in\omega}T_n$} ${\cal M}_n\models
T_n$, $n\in\omega$.}
\end{definition}

Clearly, the theory $T_1\sqcup T_2$ does not depend on choice of
disjoint union ${\cal M}_1\sqcup {\cal M}_2$ of models ${\cal
M}_1\models T_1$ and ${\cal M}_2\models T_2$. Besides, the
cardinality of ${\rm RK}(T_1\sqcup T_2)$ is equal to the product
of cardinalities for ${\rm RK}(T_1)$ and ${\rm RK}(T_2)$, and the
relation $\leq_{\rm RK}$ on ${\rm RK}(T_1\sqcup T_2)$ equals the
Pareto relation \cite{SO1} defined by preorders in  ${\rm
RK}(T_1)$ and ${\rm RK}(T_2)$. Indeed, each type $p(\bar{x})$ of
$T_1\sqcup T_2$ is isolated by set consisting of some types
$p_1(\bar{x}^1)$ and $p_2(\bar{x}^2)$ of theories $T_1$ and $T_2$
respectively, as well as of formulas $P^1(x^1_i)$ and $P^2(x^2_j)$
for all coordinates in tuples $\bar{x}^1$ and $\bar{x}^2$. For
types $p(\bar{x})$ and $p'(\bar{y})$ of $T_1\sqcup T_2$, we have
$p(\bar{x})\leq_{\rm RK} p'(\bar{y})$ if and only if
$p_1(\bar{x}^1)\leq_{\rm RK} p'_1(\bar{y}^1)$ (in $T_1$) and
$p_2(\bar{x}^2)\leq_{\rm RK} p'_2(\bar{y}^2)$ (in $T_2$).

Thus, the following proposition holds.

\medskip
\begin{proposition}\label{st676} {\rm \cite{SuCCMCT2, Su09Irk}}
For any small theories $T_1$ and $T_2$ of disjoint predicate
languages $\Sigma_1$ and $\Sigma_2$ respectively, the theory
$T_1\sqcup T_2$ is mutually ${\rm RK}$-coordinated with respect to
its restrictions to $\Sigma_1$ and $\Sigma_2$. The cardinality of
${\rm RK}(T_1\sqcup T_2)$ is equal to the product of cardinalities
for ${\rm RK}(T_1)$ and ${\rm RK}(T_2)$, i.~e.,
\begin{equation}\label{s64} I_p(T_1\sqcup
T_2,\omega)=I_p(T_1,\omega)\cdot I_p(T_2,\omega),
\end{equation} and the relation
$\leq_{\rm RK}$ on ${\rm RK}(T_1\sqcup T_2)$ equals the Pareto
relation defined by preorders in ${\rm RK}(T_1)$ and ${\rm
RK}(T_2)$.
\end{proposition}

\begin{note}\label{no6742}
{\rm \cite{SuCCMCT2, Su09Irk} An isomorphism \ of limit \ models \
of theory \ $T_1\sqcup T_2$ is defined by isomorphisms of
restrictions of these models to the sets $P_1$ and $P_2$. In this
case, a countable model is limit if and only if some its
restriction (to $P_1$ or to $P_2$) is limit and the following
equality holds:
\begin{equation}\label{s65}
I(T_1\sqcup T_2,\omega)=I(T_1,\omega)\cdot I(T_2,\omega).
\end{equation}
Thus, the operation $\sqcup$ preserves both $p$-Ehrenfeuchtness
and $l$-Ehrenfeuchtness (if components are $p$-Ehrenfeucht), and,
by $(\ref{s65})$, we obtain the equality
\begin{equation}\label{s66}
I_l(T_1\sqcup T_2)=I_l(T_1)\cdot I_p(T_2,\omega)+
I_p(T_1,\omega)\cdot I_l(T_2)+I_l(T_1)\cdot I_l(T_2).
\end{equation}}
\end{note}

\section{$o$-minimal and quite $o$-minimal theories}
\noindent

Recall \cite{PilSt} that a linearly ordered structure ${\cal M}$
is {\em $o$-minimal}\index{Structure!$o$-minimal} if each formula
definable subset of $M$ is a finite union of singletons and open
intervals $(a,b)$, where $a\in M\cup\{-\infty\}$, $b\in
M\cup\{+\infty\}$. A theory $T$ is {\em
$o$-minimal}\index{Theory!$o$-minimal} if each model of $T$ is
$o$-minimal.

As examples of Ehrenfeucht $o$-minimal theories, we mention the
theories $T^1\rightleftharpoons{\rm Th}(({\mathbb
Q};<,c_n)_{n\in\omega}$ and $T^2\rightleftharpoons{\rm
Th}(({\mathbb Q};<,c_n,c'_n)_{n\in\omega}$, where $<$ is an
ordinary strict order on the set ${\mathbb Q}$ of rationals,
constants $c_n$ form a strictly increasing sequence, and constants
$c'_n$ form a strictly decreasing sequence, $c_n<c'_n$,
$n\in\omega$.

The theory $T^1$ is an Ehrenfeucht's example \cite{Va} with
$I(T^1,\omega)=3$. It has two almost prime models and one limit
model:

{\small $\bullet$} a prime model with empty set of realizations of
type $p(x)$ isolated by the set $\{c_n<x\mid n\in\omega\}$ of
formulas;

{\small $\bullet$} a prime model over a realization of the type
$p(x)$, with the least realization of that type;

{\small $\bullet$} one limit model over the type $p(x)$, with the
set of realizations of $p(x)$ forming an open interval.

The Hasse diagram for the Rudin--Keisler preorder $\leq_{\rm RK}$
and values of the function ${\rm IL}$ of distributions of numbers
of limit models for $\sim_{\rm RK}$-classes of $T^1$ is
represented in Fig.~\ref{fig1}.

\begin{figure}[t]
\begin{center}
\unitlength 4mm
\begin{picture}(5,6.5)(-5.5,1.0)
{\footnotesize\put(2,2.5){\line(0,5){5}}
\put(2,2.5){\makebox(0,0)[cc]{$\bullet$}}
\put(2,7.5){\makebox(0,0)[cc]{$\bullet$}}
\put(2,2.5){\circle{0.6}} \put(2,7.5){\circle{0.6}}
\put(3,2.5){\makebox(0,0)[cr]{$0$}}
\put(3,7.5){\makebox(0,0)[cr]{$1$}} }
\end{picture}
\hfill \unitlength 5mm
\begin{picture}(22,7)(4,1.7)
{\footnotesize \put(19.5,2.5){\line(0,5){5}}
\put(19.5,2.5){\makebox(0,0)[cc]{$\bullet$}}
\put(19.5,5){\makebox(0,0)[cc]{$\bullet$}}
\put(19.5,7.5){\makebox(0,0)[cc]{$\bullet$}}
\put(19.5,2.5){\circle{0.6}} \put(19.5,5){\circle{0.6}}
\put(19.5,7.5){\circle{0.6}}
\put(20.5,2.5){\makebox(0,0)[cr]{$0$}}
\put(20.5,5){\makebox(0,0)[cr]{$0$}}
\put(20.5,7.5){\makebox(0,0)[cr]{$3$}} }
\end{picture}
\end{center}
\parbox[t]{0.4\textwidth}{\caption{}\label{fig1}}
\hfill
\parbox[t]{0.4\textwidth}{\caption{}\label{fig2}}

\end{figure}

The theory $T^2$ has six pairwise non-isomorphic countable models:

{\small $\bullet$} a prime model with empty set of realizations of
type $p(x)$ isolated by the set $\{c_n<x\mid
n\in\omega\}\cup\{x<c'_n\mid n\in\omega\}$;

{\small $\bullet$} a prime model over a realization of $p(x)$, with
a unique realization of this type;

{\small $\bullet$} a prime model over a realization of type
$q(x,y)$ isolated by the set $p(x)\cup p(y)\cup\{x<y\}$; here the
set of realizations of $q(x,y)$ forms a closed interval $[a,b]$;

{\small $\bullet$} three limit models over the type $q(x,y)$, in
which the sets of realizations of $q(x,y)$ are intervals of forms
$(a,b]$, $[a,b)$, $(a,b)$ respectively.

In Figure \ref{fig2} we represent the Hasse diagram
of~Rudin--Keisler preorders $\leq_{\rm RK}$ and values of
distribution functions ${\rm IL}$ of numbers of limit models on
$\sim_{\rm RK}$-equivalence classes for the theory $T^2$.

The following theorem shows that the number of countable models of
Ehrenfeucht $o$-minimal theories is exhausted by combinations of
these numbers for the theories $T_1$ and $T_2$.

\begin{theorem}\label{th1} {\rm \cite{Ma}}
Let $T$ be an $o$-minimal theory in a countable language. Then
either $T$ has $2^{\omega}$ countable models or $T$ has exactly
$3^k\cdot 6^s$ countable models, where $k$ and $s$ are natural
numbers. Moreover, for any $k,s\in\omega$ there is an $o$-minimal
theory $T$ with exactly $3^k\cdot 6^s$ countable models.
\end{theorem}

The notion of {\it weak o-minimality} was initially deeply studied
by D.~Mac\-pher\-son, D.~Mar\-ker, and C.~Steinhorn in \cite{mms}. A
subset $A$ of a linearly ordered structure $M$ is {\it convex} if
for any $a, b\in A$ and $c\in M$ whenever $a<c<b$ we have $c\in
A$. A {\it weakly o-minimal structure} is a linearly ordered
structure $\mathcal{M}=\langle M,=,<,\ldots \rangle$ such that any
definable (with parameters) subset of the structure $\mathcal{M}$
is a finite union of convex sets in $\mathcal{M}$. Real closed
fields with a proper convex valuation ring provide an important
example of weakly o-minimal (not o-minimal) structures.

In the following definitions we assume that $M$ is a weakly
o-minimal structure, $A, B\subseteq M$, $M$ is $|A|^+$-saturated,
and $p, q\in S_1(A)$ are non-algebraic types.
\begin{definition}\rm (B.S. Baizhanov, \cite{bbs1})
We say that $p$ is not {\it weakly orthogonal} to $q$
($p\not\perp^w  q$) if there are an $A$-definable formula
$H(x,y)$, $a \in p(M)$, and $b_1, b_2 \in q(M)$ such that $b_1 \in
H(M,a)$ and $b_2 \not\in H(M,a)$.
\end{definition}

\begin{lemma}\label{leq} {\rm (\cite{bbs1},
Corollary 34 (iii))}   The relation $\not\perp^w$ of the weak
non-orthogonality is an equivalence relation on $S_1(A)$.
\end{lemma}

In \cite{k2003}, quite o-minimal theories were introduced forming
a subclass of the class of weakly o-minimal theories and
preserving a series of properties for o-minimal theories. For
instance, in \cite{jms}, $\aleph_0$-categorical quite o-minimal
theories were completely described. This description implies their
binarity (the similar result holds for $\aleph_0$-categorical
o-minimal theories).
\begin{definition} \rm \cite{k2003} We say that $p$ is not
{\it quite orthogonal} to $q$ ($p\not\perp^q q$) if there is an
$A$-definable bijection $f: p(M)\to q(M)$. We say that a weakly
o-minimal theory is {\it quite o-minimal} if the relations of weak
and quite orthogonality coincide for 1-types over arbitrary sets
of models of the given theory.
\end{definition}

Clearly, any $o$-minimal theory is quite $o$-minimal, since for
non-weakly orthogonal $1$-types over an arbitrary set $A$ there is
an $A$-definable strictly monotone bijection between sets of
realizations of these types.

\begin{example}\rm  Let $\mathcal{M}=\langle M, <,P^1_1, P^1_2, E^2_1,
E^2_2, f^1\rangle$ be a linearly ordered structure such that
$\mathcal{M}$ is a disjoint union of interpretations of unary
predicates $P_1$ and $P_2$, where
$P_1(\mathcal{M})<P_2(\mathcal{M})$. We identify the
interpretations of $P_1$ and $P_2$ with $\mathbb Q\times\mathbb Q$
having the lexicographical order. For the interpretations of
binary predicates $E_1(x,y)$ and $E_2(x,y)$ we take equivalence
relations on $P_1(\mathcal{M})$ and $P_2(\mathcal{M})$,
respectively, such that for every $x=(n_1,m_1),
y=(n_2,m_2)\in\mathbb{Q}\times\mathbb{Q}$,
$$E_i(x,y)\Leftrightarrow n_1=n_2, \mbox{ where }
i=1,2.$$

The symbol $f$ is interpreted by partial unary function with
$\Dom(f)=P_1(\mathcal{M})$ and $\Range(f)=P_2(\mathcal{M})$ such
that $f((n,m))=(n,-m)$ for all $(n,m)\in
\mathbb{Q}\times\mathbb{Q}$.

It is easy to see that $E_1(x,y)$ and $E_2(x,y)$ are
$\emptyset$-definable equivalence relations dividing
$P_1(\mathcal{M})$ and $P_2(\mathcal{M})$, respectively, into
infinitely many infinite convex classes. We assert that $f$ is
strictly decreasing on each class $E_1(a,\mathcal{M})$, where
$a\in P_1(\mathcal{M})$, and  $f$ is strictly increasing on
$P_1(\mathcal{M})/E_1$. It is clear that ${\rm Th}(\mathcal{M})$
is a quite $o$-minimal theory. The theory ${\rm Th}(\mathcal{M})$
is not $o$-minimal, since $E_1(a,M)$ defines a convex set which is
not a union of finitely many intervals in $\mathcal{M}$.
\end{example}

The following theorem, proved in \cite{KS}, strengthens Theorem
\ref{th1}.

\begin{theorem}\label{th2}
Let $T$ be a quite $o$-minimal theory in a countable language.
Then either $T$ has $2^{\omega}$ countable models or $T$ has
exactly $3^k\cdot 6^s$ countable models, where $k$ and $s$ are
natural numbers. Moreover, for any $k,s\in\omega$ there is an
$o$-minimal theory $T$ with exactly $3^k\cdot 6^s$ countable
models.
\end{theorem}

It was shown in \cite{KS} that quite $o$-minimal Ehrenfeucht
theories are binary. But it does not hold in general:

\begin{example}\rm Let $\mathcal{M}=\langle M; <, P^1_1, P^1_2, P^1_3,
f^2\rangle$ be a linearly ordered structure such that $M$ is a
disjoint union of interpretations of unary predicates $P_1, P_2$,
and $P_3$, where
$P_1(\mathcal{M})<P_2(\mathcal{M})<P_3(\mathcal{M})$. We identify
each interpretation of $P_i$ ($1\le i\le 3$) with the set
$\mathbb{Q}$ of rational numbers, with ordinary orders. The symbol
$f$ is interpreted by partial binary function with
$\Dom(f)=P_1(\mathcal{M})\times P_2(\mathcal{M})$ and
$\Range(f)=P_3(\mathcal{M})$ such that $f(a,b)=a+b$ for all
$(a,b)\in\mathbb{Q}\times\mathbb{Q}$.

Clearly, ${\rm Th}(\mathcal{M})$ is a quite $o$-minimal theory.
Take arbitrary $a\in P_1(\mathcal{M}), b\in P_2(\mathcal{M})$.
Obviously, the functions $f_b(x):=f(x,b)$ and $g_a(y):=f(a,y)$ are
strictly increasing on $P_1(\mathcal{M})$ and $P_2(\mathcal{M})$,
respectively. Take an arbitrary $a_1\in P_1(\mathcal{M})$ with
$a<a_1$ and consider the following formulas:
$$\Phi_1(y,a,a_1,b):=f_b(a)=f_y(a_1)\land P_2(y),$$
$$\Phi_n(y,a,a_1,b):=\exists y_0 [\Phi_{n-1}(y_0, a,a_1,b)\land f_{y_0}(a)=f_y(a_1)\land P_2(y)],\quad n\ge 2.$$
Clearly, $\mathcal{M}\models \exists ! y \Phi_n(y,a,a_1,b)$ for
each $n<\omega$, i.e., $\dcl(\{a,a_1,b\})$ infinite. Then
considering the following set of formulas:
$$\{P_2(x)\}\cup\{x<b\}\cup \{\forall y[\Phi_n(y,a,a_1,b)\to x<y]\mid n\in \omega\}$$
and checking its local consistency, we obtain that there exists a
non-principal 1-type over $\{a,a_1,b\}$ extending the given set of
formulas. Whence, ${\rm Th}(\mathcal{M})$ has $2^{\omega}$
countable models. Since for each finite set $A\subseteq M$ there
are only at most countably many 1-types over $A$, we conclude that
the theory ${\rm Th}\mathcal{M})$ is small.
\end{example}

Thus, the following proposition is proved:

\begin{proposition}
There exists a small quite $o$-minimal theory, which is not
binary.
\end{proposition}

\begin{definition}
{\rm \cite{SuCCMCT1, Su09Irk} We say that small theories $T_1$ and
$T_2$ are {\em characteristically
equivalent}\index{Theories!characteristically equivalent} and
write $T_1\sim_{\rm ch}T_2$\index{$T_1\sim_{\rm ch}T_2$} if the
structure ${\rm RK}(T_1)$ is isomorphic to the structure ${\rm
RK}(T_2)$ and, by the corresponding replacement of isomorphism
types in ${\rm RK}(T_1)$ to isomorphism types in ${\rm RK}(T_2)$,
the distribution function {\rm IL} for numbers of limit models of
$T_1$ is transformed to the distribution function for numbers of
limit models of $T_2$.}
\end{definition}

The following theorem is a reformulation of Theorem \ref{th2} for
quite $o$-minimal Ehrenfeucht theories producing the direct
generalization of Theorem 1.1.5.3 in \cite{SuCCMCT1}.

\medskip
\begin{theorem}\label{thchar}
Any model of a quite $o$-minimal Ehrenfeucht theory $T$ is densely
ordered besides, possibly, finitely many elements with successors
or predecessors laying in the definable closure of empty set. The
theory $T$ is characteristically equivalent to some finite
disjoint union of theories of form $T^1$, $T^2$ {\rm (}$T\sim_{\rm
ch}\bigsqcup\limits_{i=1}^k T^1_i\sqcup\bigsqcup\limits_{j=1}^l
T^2_j$, where $T^1_i$ are similar to $T^1$ and $T^2_j$ are similar
to $T^2${\rm )} and has $3^k\cdot 6^l$ pairwise non-isomorphic
countable models.
\end{theorem}

\section{Distributions of countable models}
\noindent

In this section, using Theorems \ref{th1121_1136} and \ref{thchar}
we give a description of Rudin--Keisler preorders and distribution
functions of numbers of limit models for quite $o$-minimal
Ehrenfeucht theories, as well as propose representations of this
distributions, based on the decomposition formula (\ref{eqmain}).

In view of Proposition \ref{st676} and Theorem \ref{thchar} the
Hasse diagrams for distributions of countable models for quite
$o$-minimal Ehrenfeucht theories are constructed as figures of
Pareto relations for disjoint unions of copies of theories $T^1$
and $T^2$, i.e., they are combinations of the Hasse diagrams shown
in Fig.~\ref{fig1} and \ref{fig2}.

Now we describe the distributions above for the theories
$\bigsqcup\limits_{i=1}^k T^1_i$.

In Fig.~\ref{fig3} and \ref{fig4} the Hasse diagrams are shown for
the theories $T^1_1\sqcup T^1_2$ and $T^1_1\sqcup T^1_2\sqcup
T^1_3$, respectively.

\begin{figure}[t]
\begin{center}
\unitlength 14mm
\begin{picture}(5,1)(-1.14,-0.1)
{\footnotesize \put(0,1){\makebox(0,0)[cc]{$\bullet$}}
\put(0,1){\circle{0.2}} \put(1,0){\circle{0.2}}
\put(1,0){\makebox(0,0)[cc]{$\bullet$}}
\put(2,1){\makebox(0,0)[cc]{$\bullet$}}
\put(2,1){\circle{0.2}}\put(1,2){\circle{0.2}}
\put(1,2){\makebox(0,0)[cc]{$\bullet$}} \put(1,0){\line(1,1){1}}
\put(1,0){\line(-1,1){1}} \put(1,2){\line(-1,-1){1}}
\put(1,2){\line(1,-1){1}} \put(1,-0.2){\makebox(0,0)[tc]{0}}
\put(1,2.2){\makebox(0,0)[bc]{$3$}}
\put(0.06,1.3){\makebox(0,0)[cr]{$1$}}
\put(1.95,1.3){\makebox(0,0)[cl]{$1$}} }
\end{picture}
\hfill \unitlength 14mm
\begin{picture}(5,3.5)(1.35,0.2)
{\footnotesize\put(3,1){\makebox(0,0)[cc]{$\bullet$}}\put(3,1){\circle{0.2}}
\put(3,2){\makebox(0,0)[cc]{$\bullet$}}\put(3,2){\circle{0.2}}
\put(4,0){\makebox(0,0)[cc]{$\bullet$}}\put(4,0){\circle{0.2}}
\put(4,1){\makebox(0,0)[cc]{$\bullet$}}\put(4,1){\circle{0.2}}
\put(4,2){\makebox(0,0)[cc]{$\bullet$}}\put(4,2){\circle{0.2}}
\put(4,3){\makebox(0,0)[cc]{$\bullet$}}\put(4,3){\circle{0.2}}
\put(5,1){\makebox(0,0)[cc]{$\bullet$}}\put(5,1){\circle{0.2}}
\put(5,2){\makebox(0,0)[cc]{$\bullet$}}\put(5,2){\circle{0.2}}
\put(3,1){\line(0,1){1}} \put(5,1){\line(0,1){1}}
\put(4,0){\line(0,1){1}} \put(4,2){\line(0,1){1}}
\put(4,0){\line(1,1){1}} \put(4,1){\line(1,1){1}}
\put(4,0){\line(-1,1){1}} \put(4,1){\line(-1,1){1}}
\put(4,3){\line(-1,-1){1}} \put(4,2){\line(-1,-1){1}}
\put(4,3){\line(1,-1){1}} \put(4,2){\line(1,-1){1}}
\put(4.2,0){\makebox(0,0)[cl]{0}}
\put(3.2,1){\makebox(0,0)[cl]{1}}
\put(4.2,1){\makebox(0,0)[cl]{1}}
\put(5.2,1){\makebox(0,0)[cl]{1}}
\put(3.2,2){\makebox(0,0)[cl]{3}}
\put(4.2,2){\makebox(0,0)[cl]{3}}
\put(5.2,2){\makebox(0,0)[cl]{3}}
\put(4.2,3.1){\makebox(0,0)[cl]{7}}}
\end{picture}
\end{center}
\parbox[t]{0.4\textwidth}{\caption{}\label{fig3}}
\hfill
\parbox[t]{0.4\textwidth}{\caption{}\label{fig4}}

\end{figure}

Adding new disjoint copies of $T^1$ we note that ${\rm RK}(T)$,
where $T=\bigsqcup\limits_{i=1}^k T^1_i$, forms a $k$-dimensional
cube $Q_k$ \cite{SO1}, i.å., represented as a finite Boolean
algebra $\mathcal{B}_k$ with $k$ atoms $u_1,\ldots,u_k$. These
atoms correspond to models realizing unique $1$-types in the set
$\{p_1(x)$, $\ldots$, $p_k(x)\}$ of all nonprincipal $1$-types. Thus,
each element $u_{i_1}\vee\ldots\vee u_{i_t}$ of the Boolean
algebra $\mathcal{B}_k$ corresponds to an almost prime model of
$T$, realizing only nonprincipal $1$-types
$p_{i_1}(x),\ldots,p_{i_t}(x)$.

The number of limit models for the element $u_{i_1}\vee\ldots\vee
u_{i_t}$, i.~e., of limit models over (unique) completion
$q_{i_1,\ldots,i_t}(x_1,\ldots,x_t)$ of the type
$p_{i_1}(x_1)\cup\ldots\cup p_{i_t}(x_t)$ equals $2^t-1$. Indeed,
choosing a prime model over the type $q_{i_1,\ldots,i_t}$ we have
$2^t$ possibilities characterizing an independent choice either
prime or limit model over each type $p_{i_j}$. Removing the
(unique) possibility of choice of prime model for each type
$p_{i_j}$, i.~e., of prime model over the type
$q_{i_1,\ldots,i_t}$, we obtain the value
\begin{equation}\label{qeq3}
I_l(T,q_{i_1,\ldots,i_t})=2^t-1
\end{equation}
of the number of limit models over the type $q_{i_1,\ldots,i_t}$.

Since there are $3^k$ countable models, $2^k$ of them are almost
prime, and the remaining are limit ones, the total number of limit
models, calculated on the basis of relations (\ref{qeq3}) (see
also (\ref{s66})) leads to the following:
\begin{equation}\label{qeq4}
\sum\limits_{q_{i_1,\ldots,i_m}}
I_l(T,q_{i_1,\ldots,i_t})=\sum\limits_{t=1}^k(2^t-1)\cdot
C^t_k=3^k-2^k.
\end{equation}
By (\ref{qeq4}) for the theories $\bigsqcup\limits_{i=1}^k T^1_i$,
we have the following representation of the decomposition formula
(\ref{eqmain}):
\begin{equation}\label{qeq5}
3^k=2^k+\sum\limits_{t=1}^k(2^t-1)\cdot C^t_k.
\end{equation}

For $k=1$ we have $3=2+1$, for $k=2$: \ $9=4+1\cdot 2+3\cdot 1$,
for $k=3$: \ $27=8+1\cdot 3+3\cdot 3+7\cdot 1$, as shown in Fig.
\ref{fig1}, \ref{fig3}, \ref{fig4}, respectively.

Notice that each $(k+1)$-th diagram contains $k+1$ previous ones.

\medskip
Now we describe the distributions for the theories
$\bigsqcup\limits_{j=1}^s T^2_j$.

In Fig.~\ref{fig5} and \ref{fig6} the Hasse diagrams shown for the
theories $T^2_1\sqcup T^2_2$ and $T^2_1\sqcup T^2_2\sqcup T^2_3$,
respectively. These diagrams form non-distributive lattices, which
are obtained, respectively, from four-element and eight-element
Boolean algebras by extensions of each two-dimensional cube by
four new elements such that each edge of given Boolean algebra
contains new intermediate element. The theory $T^2_1\sqcup T^2_2$
has $6^2=36$ countable models, where $9$ of them are almost prime
and $27$ are limit. The theory $T^2_1\sqcup T^2_2\sqcup T^2_3$ has
$6^3=216$ countable models, where $27$ of them are almost prime
and $189$ are limit.

\begin{figure}
\begin{center}
\unitlength 18mm
\begin{picture}(3,1)(-0.69,-0.1)
{\footnotesize \put(0,1){\makebox(0,0)[cc]{$\bullet$}}
\put(0,1){\circle{0.16}} \put(1,0){\circle{0.16}}
\put(1,0){\makebox(0,0)[cc]{$\bullet$}}
\put(2,1){\makebox(0,0)[cc]{$\bullet$}}
\put(2,1){\circle{0.16}}\put(1,2){\circle{0.16}}
\put(1,2){\makebox(0,0)[cc]{$\bullet$}}
\put(1,0){\line(1,1){1}}\put(1.5,0.5){\line(-1,1){1}}
\put(0.5,0.5){\line(1,1){1}} \put(1,0){\line(-1,1){1}}
\put(1,2){\line(-1,-1){1}} \put(1,2){\line(1,-1){1}}
\put(1.5,0.5){\makebox(0,0)[cc]{$\bullet$}}
\put(0.5,0.5){\makebox(0,0)[cc]{$\bullet$}}
\put(1,1){\makebox(0,0)[cc]{$\bullet$}}
\put(0.5,1.5){\makebox(0,0)[cc]{$\bullet$}}
\put(1.5,1.5){\makebox(0,0)[cc]{$\bullet$}}
\put(1.5,0.5){\circle{0.16}} \put(0.5,0.5){\circle{0.16}}
\put(1,1){\circle{0.16}} \put(0.5,1.5){\circle{0.16}}
\put(1.5,1.5){\circle{0.16}} \put(1,-0.17){\makebox(0,0)[tc]{0}}
\put(1,2.15){\makebox(0,0)[bc]{$15$}}
\put(0.06,1.22){\makebox(0,0)[cr]{$3$}}
\put(0.56,0.7){\makebox(0,0)[cr]{$0$}}
\put(1.56,0.7){\makebox(0,0)[cr]{$0$}}
\put(1.56,1.7){\makebox(0,0)[cr]{$3$}}
\put(0.56,1.7){\makebox(0,0)[cr]{$3$}}
\put(1.06,1.22){\makebox(0,0)[cr]{$0$}}
\put(1.95,1.22){\makebox(0,0)[cl]{$3$}} }
\end{picture}
\hfill \unitlength 20mm
\begin{picture}(4,3.5)(1.6,0.2)
{\footnotesize\put(3,1){\makebox(0,0)[cc]{$\bullet$}}\put(3,1){\circle{0.14}}
\put(3,2){\makebox(0,0)[cc]{$\bullet$}}\put(3,2){\circle{0.14}}
\put(4,0){\makebox(0,0)[cc]{$\bullet$}}\put(4,0){\circle{0.14}}
\put(4,1){\makebox(0,0)[cc]{$\bullet$}}\put(4,1){\circle{0.14}}
\put(4.006,2){\makebox(0,0)[cc]{$\bullet$}}\put(4.006,2){\circle{0.08}}
\put(4,3){\makebox(0,0)[cc]{$\bullet$}}\put(4,3){\circle{0.14}}
\put(5,1){\makebox(0,0)[cc]{$\bullet$}}\put(5,1){\circle{0.14}}
\put(5,2){\makebox(0,0)[cc]{$\bullet$}}\put(5,2){\circle{0.14}}
\put(3,1){\line(0,1){1}} \put(5,1){\line(0,1){1}}
\put(4,0){\line(0,1){1}} \put(4,2){\line(0,1){1}}
\put(4,0){\line(1,1){1}} \put(4,1){\line(1,1){1}}
\put(4,0){\line(-1,1){1}} \put(4,1){\line(-1,1){1}}
\put(4,3){\line(-1,-1){1}} \put(4,2){\line(-1,-1){1}}
\put(4,3){\line(1,-1){1}} \put(4,2){\line(1,-1){1}}
\put(4,2.5){\makebox(0,0)[cc]{$\bullet$}}\put(4,2.5){\circle{0.14}}
\put(3.35,0.65){\makebox(0,0)[cc]{$\bullet$}}\put(3.35,0.65){\circle{0.14}}
\put(3.35,1.15){\makebox(0,0)[cc]{$\bullet$}}\put(3.35,1.15){\circle{0.14}}
\put(3.35,1.65){\makebox(0,0)[cc]{$\bullet$}}\put(3.35,1.65){\circle{0.14}}
\put(3.95,1.75){\makebox(0,0)[cc]{$\bullet$}}\put(3.95,1.75){\circle{0.14}}
\put(3.95,1.25){\makebox(0,0)[cc]{$\bullet$}}\put(3.95,1.25){\circle{0.14}}
\put(3.945,2.25){\makebox(0,0)[cc]{$\bullet$}}\put(3.945,2.25){\circle{0.08}}
\put(4.35,1.65){\makebox(0,0)[cc]{$\bullet$}}\put(4.35,1.65){\circle{0.14}}
\put(4.35,2.15){\makebox(0,0)[cc]{$\bullet$}}\put(4.35,2.15){\circle{0.14}}
\put(4.35,2.65){\makebox(0,0)[cc]{$\bullet$}}\put(4.35,2.65){\circle{0.14}}
\put(3.6,1.6){\makebox(0,0)[cc]{$\bullet$}}\put(3.6,1.6){\circle{0.14}}
\put(3.6,2.1){\makebox(0,0)[cc]{$\bullet$}}\put(3.6,2.1){\circle{0.14}}
\put(3.6,2.6){\makebox(0,0)[cc]{$\bullet$}}\put(3.6,2.6){\circle{0.14}}
\put(4.6,0.6){\makebox(0,0)[cc]{$\bullet$}}\put(4.6,0.6){\circle{0.14}}
\put(4.6,1.1){\makebox(0,0)[cc]{$\bullet$}}\put(4.6,1.1){\circle{0.14}}
\put(4.6,1.6){\makebox(0,0)[cc]{$\bullet$}}\put(4.6,1.6){\circle{0.14}}
\put(4,0.5){\makebox(0,0)[cc]{$\bullet$}}\put(4,0.5){\circle{0.14}}
\put(5,1.5){\makebox(0,0)[cc]{$\bullet$}}\put(5,1.5){\circle{0.14}}
\put(3,1.5){\makebox(0,0)[cc]{$\bullet$}}\put(3,1.5){\circle{0.14}}

\put(3.35,0.65){\line(0,1){1}} \put(4.6,0.6){\line(0,1){1}}
\put(3.35,0.65){\line(1,1){1}} \put(4.6,0.6){\line(-1,1){1}}
\put(4,0.5){\line(1,1){1}} \put(4,0.5){\line(-1,1){1}}
\put(3.35,1.65){\line(1,1){1}} \put(4.6,1.6){\line(-1,1){1}}
\put(3,1.5){\line(1,1){1}} \put(5,1.5){\line(-1,1){1}}
\put(3.6,1.6){\line(0,1){1}} \put(4.35,1.65){\line(0,1){1}}
\put(3.35,1.15){\line(1,1){1}} \put(4.6,1.1){\line(-1,1){1}}
\put(3.95,1.25){\line(0,1){1}}

\put(4.1,0){\makebox(0,0)[cl]{0}}
\put(4.1,0.5){\makebox(0,0)[cl]{0}}
\put(4.69,0.6){\makebox(0,0)[cl]{0}}
\put(4.69,1.1){\makebox(0,0)[cl]{0}}
\put(4.69,1.6){\makebox(0,0)[cl]{3}}
\put(3.17,0.65){\makebox(0,0)[cl]{0}}
\put(3.193,1.14){\makebox(0,0)[cl]{0}}
\put(3.193,1.64){\makebox(0,0)[cl]{3}}
\put(3.445,1.63){\makebox(0,0)[cl]{3}}
\put(4.185,1.66){\makebox(0,0)[cl]{3}}

\put(2.83,1){\makebox(0,0)[cl]{3}}
\put(2.83,1.5){\makebox(0,0)[cl]{3}}
\put(2.73,2){\makebox(0,0)[cl]{15}}
\put(4.076,1){\makebox(0,0)[cl]{3}}
\put(5.1,1){\makebox(0,0)[cl]{3}}
\put(5.1,1.5){\makebox(0,0)[cl]{3}}
\put(5.1,2){\makebox(0,0)[cl]{15}}
\put(4.04,1.25){\makebox(0,0)[cl]{0}}
\put(4.04,1.75){\makebox(0,0)[cl]{0}}
\put(3.82,2.0){\makebox(0,0)[cl]{{\tiny 15}}}
\put(3.44,2.1){\makebox(0,0)[cl]{3}}
\put(3.36,2.6){\makebox(0,0)[cl]{15}}
\put(3.765,2.52){\makebox(0,0)[cl]{15}}
\put(4.01,2.25){\makebox(0,0)[cl]{{\tiny 3}}}
\put(4.43,2.15){\makebox(0,0)[cl]{3}}
\put(4.43,2.65){\makebox(0,0)[cl]{15}}
\put(3.929,3.16){\makebox(0,0)[cl]{63}}

}
\end{picture}
\end{center}
\parbox[t]{0.4\textwidth}{\caption{}\label{fig5}}
\hfill
\parbox[t]{0.4\textwidth}{\caption{}\label{fig6}}

\end{figure}

Continuing the process of adding disjoint copies of the theory
$T^2$, we observe that  ${\rm RK}(T)$, where
$T=\bigsqcup\limits_{j=1}^s T^2_j$, is obtained from
$s$-dimensional cube replacing edges by three-element lines and
forming $s$-dimensional linear space ${\mathcal L}_{s,3}$ over the
field $\mathbb Z_3$. Therefore, $|{\rm RK}(T)|=3^s$. Here, the
theory $T$ has exactly $s$ nonprincipal $1$-types
$p_1(x),\ldots,p_s(x)$, each of which, in almost prime models,
either does not have realizations, or has unique realization, or
has infinitely many realizations including the least and the
greatest ones.

To calculate the number of limit models, we note that the
structure ${\mathcal L}_{s,3}$ contains the $s$-dimensional cube,
whose vertices, \ $2^s$ \ ones, symbolize prime models over
completions \ $q_{j_1,\ldots,j_m}(x_1,\ldots,x_m)$ of types
$p_{j_1}(x_1)\cup\ldots\cup p_{j_m}(x_m)$ such that these prime
models have at most one realization for each type
$p_1(x),\ldots,p_s(x)$ and do not generate limit models.
Furthermore, we choose among $s$ types $p_j$ some $m$ types,
responsible for the existence of limit models generated by
realizations of these types, and obtain $4^m-1$ possibilities for
these limit models by variations of existence or absence of least
and greatest realizations. Together with the choice of $m$ types
we choose among remaining $s-m$ types some $r$ types having unique
realizations. Under these conditions of choice we have
$(4^m-1)\cdot C^m_s\cdot C^r_{s-m}$ possibilities. Summarizing
these values we obtain the following equations:
$$\sum\limits_{q_{i_1,\ldots,i_m}}
I_l(T,q_{i_1,\ldots,i_m})=\sum\limits_{m=1}^s\sum\limits_{r=0}^{s-m}
(4^m-1)\cdot C^m_s\cdot C^r_{s-m}=$$
\begin{equation}\label{qeq6}
=\sum\limits_{m=1}^s\left(\sum\limits_{r=0}^{s-m}
C^r_{s-m}\right)(4^m-1)\cdot C^m_s=\sum\limits_{m=1}^s
2^{s-m}\cdot (4^m-1)\cdot C^m_s=6^s-3^s.
\end{equation}
By (\ref{qeq6}) for the theory $\bigsqcup\limits_{j=1}^s T^2_j$,
we have the following representation of the decomposition formula
(\ref{eqmain}):
\begin{equation}\label{qeq7}
6^s=3^s+\sum\limits_{m=1}^s 2^{s-m}\cdot (4^m-1)\cdot C^m_s.
\end{equation}

For $s=1$ we have $6=3+1\cdot 3\cdot 1$, for $s=2$: \ $36=9+2\cdot
3\cdot 2 +1\cdot 15\cdot 1$, for $s=3$: \ $216=27+4\cdot 3\cdot
3+2\cdot 15\cdot 3+1\cdot 63\cdot 1$, as shown in Fig. \ref{fig2},
\ref{fig5}, \ref{fig6}, respectively.

Notice that similarly to the cases $3^k$, for the cases $6^s$ each
successive $(s + 1)$-th diagram contains $2(s+1)$ previous ones.

\medskip
Finally, we describe the indicated distributions for the theories
$\bigsqcup\limits_{i=1}^k T^1_i\sqcup\bigsqcup\limits_{j=1}^s
T^2_j$.

In Fig.~\ref{fig7}, \ref{fig8} and \ref{fig9}, the Hasse diagrams
are shown for the theories $T^1_1\sqcup T^2_1$, $T^1_1\sqcup
T^1_2\sqcup T^2_1$, and $T^1_1\sqcup T^2_1\sqcup T^2_1$,
respectively. The theory $T^1_1\sqcup T^2_1$ has $3\cdot 6=18$
countable models, $6$ of them are almost prime and $12$ are limit
ones. The theory $T^1_1\sqcup T^1_2\sqcup T^2_1$ has $3^2\cdot
6=54$ countable models, $12$ of them are almost prime and $42$ are
limit ones. The theory $T^1_1\sqcup T^2_1\sqcup T^2_1$ has $3\cdot
6^2=108$ countable models,$18$ of them are almost prime and $90$
are limit ones.

\begin{figure}[t]
\begin{center}
\unitlength 24mm
\begin{picture}(3,1)(-0.05,-0.1)
{\footnotesize  \put(1,0){\circle{0.16}}
\put(1,0){\makebox(0,0)[cc]{$\bullet$}}
\put(2,1){\makebox(0,0)[cc]{$\bullet$}} \put(2,1){\circle{0.16}}
\put(1,0){\line(1,1){1}}\put(1.5,0.5){\line(-1,1){0.5}}
\put(0.5,0.5){\line(1,1){1}} \put(1,0){\line(-1,1){0.5}}
\put(1.5,1.5){\line(1,-1){0.5}}
\put(1.5,0.5){\makebox(0,0)[cc]{$\bullet$}}
\put(0.5,0.5){\makebox(0,0)[cc]{$\bullet$}}
\put(1,1){\makebox(0,0)[cc]{$\bullet$}}

\put(1.5,1.5){\makebox(0,0)[cc]{$\bullet$}}
\put(1.5,0.5){\circle{0.16}} \put(0.5,0.5){\circle{0.16}}
\put(1,1){\circle{0.16}} \put(1.5,1.5){\circle{0.16}}
\put(1,-0.17){\makebox(0,0)[tc]{0}}

\put(0.56,0.7){\makebox(0,0)[cr]{$1$}}
\put(1.56,0.7){\makebox(0,0)[cr]{$0$}}
\put(1.56,1.7){\makebox(0,0)[cr]{$7$}}

\put(1.06,1.22){\makebox(0,0)[cr]{$1$}}
\put(1.95,1.22){\makebox(0,0)[cl]{$3$}} }
\end{picture}
\hfill \unitlength 24mm
\begin{picture}(3,2.5)(2.65,0.1)
{\footnotesize

\put(4,0){\makebox(0,0)[cc]{$\bullet$}}\put(4,0){\circle{0.14}}
\put(4,1){\makebox(0,0)[cc]{$\bullet$}}\put(4,1){\circle{0.14}}

\put(5,1){\makebox(0,0)[cc]{$\bullet$}}\put(5,1){\circle{0.14}}

\put(5,1){\line(0,1){0.5}} \put(4,0){\line(0,1){0.5}}
 \put(4,0){\line(1,1){1}}
\put(4,0){\line(-1,1){0.5}}

 \put(4.5,1.5){\line(1,-1){0.5}}

\put(3.5,0.5){\makebox(0,0)[cc]{$\bullet$}}\put(3.5,0.5){\circle{0.14}}
\put(3.5,1.0){\makebox(0,0)[cc]{$\bullet$}}\put(3.5,1.0){\circle{0.14}}

\put(4.0,1.5){\makebox(0,0)[cc]{$\bullet$}}\put(4.0,1.5){\circle{0.14}}
\put(4.5,1.5){\makebox(0,0)[cc]{$\bullet$}}\put(4.5,1.5){\circle{0.14}}
\put(4.5,2.0){\makebox(0,0)[cc]{$\bullet$}}\put(4.5,2.0){\circle{0.14}}

\put(4.5,0.5){\makebox(0,0)[cc]{$\bullet$}}\put(4.5,0.5){\circle{0.14}}
\put(4.5,1){\makebox(0,0)[cc]{$\bullet$}}\put(4.5,1){\circle{0.14}}

\put(4,0.5){\makebox(0,0)[cc]{$\bullet$}}\put(4,0.5){\circle{0.14}}
\put(5,1.5){\makebox(0,0)[cc]{$\bullet$}}\put(5,1.5){\circle{0.14}}

\put(3.5,0.5){\line(0,1){0.5}} \put(4.5,0.5){\line(0,1){0.5}}
\put(3.5,0.5){\line(1,1){1}}
\put(4.5,0.5){\line(-1,1){0.5}}\put(4.5,1){\line(-1,1){0.5}}
\put(4,0.5){\line(1,1){1}} \put(4,0.5){\line(-1,1){0.5}}

\put(5,1.5){\line(-1,1){0.5}} \put(4.5,1.5){\line(0,1){0.5}}
\put(3.5,1.0){\line(1,1){1}} \put(4,1){\line(0,1){0.5}}

\put(3.83,0){\makebox(0,0)[cl]{0}}
\put(3.83,0.5){\makebox(0,0)[cl]{1}}
\put(4.61,0.5){\makebox(0,0)[cl]{0}}
\put(4.61,1.0){\makebox(0,0)[cl]{1}}
\put(4.61,1.5){\makebox(0,0)[cl]{7}}
\put(4.45,2.17){\makebox(0,0)[cl]{15}}
\put(3.33,0.5){\makebox(0,0)[cl]{1}}
\put(3.33,1.0){\makebox(0,0)[cl]{3}}

\put(3.83,1){\makebox(0,0)[cl]{1}}
\put(3.83,1.5){\makebox(0,0)[cl]{3}}
\put(5.1,1){\makebox(0,0)[cl]{3}}
\put(5.1,1.5){\makebox(0,0)[cl]{7}}

}
\end{picture}
\end{center}
\parbox[t]{0.4\textwidth}{\caption{}\label{fig7}}
\hfill
\parbox[t]{0.4\textwidth}{\caption{}\label{fig8}}

\end{figure}

\begin{figure}[t]
\begin{center}
\unitlength 21mm
\begin{picture}(4,2.95)(2.05,0.0)
{\footnotesize\put(3,1){\makebox(0,0)[cc]{$\bullet$}}\put(3,1){\circle{0.14}}

\put(4,0){\makebox(0,0)[cc]{$\bullet$}}\put(4,0){\circle{0.14}}
\put(4,1){\makebox(0,0)[cc]{$\bullet$}}\put(4,1){\circle{0.14}}
\put(4.0,2){\makebox(0,0)[cc]{$\bullet$}}\put(4.0,2){\circle{0.14}}

\put(5,1){\makebox(0,0)[cc]{$\bullet$}}\put(5,1){\circle{0.14}}

\put(3,1){\line(0,1){0.5}} \put(5,1){\line(0,1){0.5}}
\put(4,0){\line(0,1){0.5}} \put(4,2){\line(0,1){0.5}}
\put(4,0){\line(1,1){1}} \put(4,0){\line(-1,1){1}}
 \put(4,2){\line(-1,-1){1}}
 \put(4,2){\line(1,-1){1}}
\put(4,2.5){\makebox(0,0)[cc]{$\bullet$}}\put(4,2.5){\circle{0.14}}
\put(3.5,0.5){\makebox(0,0)[cc]{$\bullet$}}\put(3.5,0.5){\circle{0.14}}
\put(3.5,1.0){\makebox(0,0)[cc]{$\bullet$}}\put(3.5,1.0){\circle{0.14}}
\put(3.5,1.5){\makebox(0,0)[cc]{$\bullet$}}\put(3.5,1.5){\circle{0.14}}
\put(3.5,2.0){\makebox(0,0)[cc]{$\bullet$}}\put(3.5,2.0){\circle{0.14}}
\put(4.0,1.5){\makebox(0,0)[cc]{$\bullet$}}\put(4.0,1.5){\circle{0.14}}
\put(4.5,1.5){\makebox(0,0)[cc]{$\bullet$}}\put(4.5,1.5){\circle{0.14}}
\put(4.5,2.0){\makebox(0,0)[cc]{$\bullet$}}\put(4.5,2.0){\circle{0.14}}

\put(4.5,0.5){\makebox(0,0)[cc]{$\bullet$}}\put(4.5,0.5){\circle{0.14}}
\put(4.5,1){\makebox(0,0)[cc]{$\bullet$}}\put(4.5,1){\circle{0.14}}

\put(4,0.5){\makebox(0,0)[cc]{$\bullet$}}\put(4,0.5){\circle{0.14}}
\put(5,1.5){\makebox(0,0)[cc]{$\bullet$}}\put(5,1.5){\circle{0.14}}
\put(3,1.5){\makebox(0,0)[cc]{$\bullet$}}\put(3,1.5){\circle{0.14}}

\put(3.5,0.5){\line(0,1){0.5}} \put(4.5,0.5){\line(0,1){0.5}}
\put(3.5,0.5){\line(1,1){1}}
\put(4.5,0.5){\line(-1,1){1}}\put(4.5,1){\line(-1,1){1}}
\put(4,0.5){\line(1,1){1}} \put(4,0.5){\line(-1,1){1}}
 \put(3,1.5){\line(1,1){1}}
\put(5,1.5){\line(-1,1){1}} \put(3.5,1.5){\line(0,1){0.5}}
\put(4.5,1.5){\line(0,1){0.5}} \put(3.5,1.0){\line(1,1){1}}
\put(4,1){\line(0,1){0.5}}

\put(4.1,0){\makebox(0,0)[cl]{0}}
\put(4.1,0.5){\makebox(0,0)[cl]{1}}
\put(4.61,0.5){\makebox(0,0)[cl]{0}}
\put(4.61,1.0){\makebox(0,0)[cl]{1}}
\put(4.61,1.5){\makebox(0,0)[cl]{3}}
\put(4.61,2.0){\makebox(0,0)[cl]{7}}
\put(3.33,0.5){\makebox(0,0)[cl]{0}}
\put(3.33,1.0){\makebox(0,0)[cl]{1}}
\put(3.33,1.5){\makebox(0,0)[cl]{3}}
\put(3.33,2.0){\makebox(0,0)[cl]{7}}

\put(2.83,1){\makebox(0,0)[cl]{3}}
\put(2.83,1.5){\makebox(0,0)[cl]{7}}

\put(4.1,1){\makebox(0,0)[cl]{0}}
\put(4.1,1.5){\makebox(0,0)[cl]{1}}
\put(5.1,1){\makebox(0,0)[cl]{3}}
\put(5.1,1.5){\makebox(0,0)[cl]{7}}

\put(4.1,2.0){\makebox(0,0)[cl]{15}}

\put(3.93,2.67){\makebox(0,0)[cl]{31}}

}
\end{picture}
\end{center}
{ \caption{} \label{fig9}}

\end{figure}

To calculate the number of limit models, we note that in the
structure ${\rm RK}(T)$, where $T=\bigsqcup\limits_{i=1}^k
T^1_i\sqcup\bigsqcup\limits_{j=1}^s T^2_j$, has the
$k$-dimensional cube $Q_k$ and the graph structure $L_{s,3}$
defined by the space ${\mathcal L}_{s,3}$. Here, the structure
${\rm RK}(T)$ is represented as the lattice with the Hasse diagram
defined by the product $Q_k\times L_{s,3}$ of graphs, and
therefore it has $2^k\cdot 3^s$ elements. Below we will also
denote the correspondent lattices by $Q_k\times L_{s,3}$.

Each vertex in  ${\rm RK}(T)$  symbolizes a prime model over
(unique) completion
$$q_{i_1,\ldots,i_t,j_1,\ldots,j_m}(x_1,\ldots,x_r,y_1,\ldots,y_m)$$
of type $p_{i_1}(x_1)\cup\ldots\cup p_{i_t}(x_m)\cup
p'_{j_1}(y_1)\cup\ldots\cup p'_{j_m}(y_m)$, where the types
$p_{1}(x),\ldots, p_{k}(x)$ exhaust the list of nonprincipal
$1$-types of the theories $T^1_i$, and the types
$p'_{1}(x),\ldots, p'_{s}(x)$ for the list of nonprincipal
$1$-types of theories $T^2_j$.  Here, almost prime models,
realizing the types $p_{i_1}(x_1),\ldots,p_{i_t}(x_m)$, have their
least realizations, as well as they have either not more than one
realizations of each type $p'_1(x),\ldots,p'_s(x)$, or, in the
latter case $p'_j(x)$, these realizations, for a fixed type, form
closed intervals.

Further, we choose among $k$ types $p_i$ some $t$ types, and among
$s$ types $p'_j$ some $m$ types, responsible for the existence of
limit models generated by realizations of these types, and obtain
$(2^t\cdot 4^m-1)$ possibilities for limit models. Together with
the choice of $m$ types we choose among remaining $s-m$ types
$p'_j$ some $r$ types having unique realizations. Under these
conditions of choice we have $(2^t\cdot 4^m-1)\cdot C^t_k\cdot
C^m_s\cdot C^r_{s-m}$ possibilities. Summarizing these values we
obtain the following equations:
$$\sum\limits_{q_{i_1,\ldots,i_t,j_1,\ldots,j_m}}
I_l(T,q_{i_1,\ldots,i_m})=\sum\limits_{t=0}^k\sum\limits_{m=0}^s\sum\limits_{r=0}^{s-m}
(2^t\cdot 4^m-1)\cdot C^t_k\cdot C^m_s\cdot C^r_{s-m}=$$
$$=\sum\limits_{t=0}^k\sum\limits_{m=0}^s\left(\sum\limits_{r=0}^{s-m}
C^r_{s-m}\right)(2^t\cdot 4^m-1)\cdot C^t_k\cdot C^m_s=$$

\begin{equation}\label{qeq8}
=\sum\limits_{t=0}^k\sum\limits_{m=0}^s 2^{s-m}\cdot(2^t\cdot
4^m-1)\cdot C^t_k\cdot C^m_s=3^k\cdot 6^s-2^k\cdot 3^s.
\end{equation}
By (\ref{qeq8}) for the theory $\bigsqcup\limits_{i=1}^k
T^1_i\sqcup\bigsqcup\limits_{j=1}^s T^2_j$, we have the following
representation of the decomposition formula (\ref{eqmain}):
\begin{equation}\label{qeq9}
3^k\cdot 6^s=2^k\cdot 3^s+\sum\limits_{t=0}^k\sum\limits_{m=0}^s
2^{s-m}\cdot(2^t\cdot 4^m-1)\cdot C^t_k\cdot C^m_s.
\end{equation}

For $k=1$ and $s=1$ we have $18=6+2\cdot 1\cdot 1\cdot 1+1\cdot
3\cdot 1\cdot 1+1\cdot 7\cdot 1\cdot 1$; for $k=2$ and $s=1$: \
$54=12+2\cdot 1\cdot 2\cdot 1 +2\cdot 3\cdot 1\cdot 1+1\cdot
3\cdot 1\cdot 1+1\cdot 7\cdot 2\cdot 1+1\cdot 15\cdot 1\cdot 1$;
for $k=1$ and $s=2$: \ $108=18+4\cdot 1\cdot 1\cdot 1+2\cdot
3\cdot 1\cdot 2+2\cdot 7\cdot 1\cdot 2+ 1\cdot 15\cdot 1\cdot
1+1\cdot 31\cdot 1\cdot 1$, as shown in Fig. \ref{fig7},
\ref{fig8}, \ref{fig9}, respectively.

By Theorem \ref{thchar} and obtained decomposition formulas
(\ref{qeq5}), (\ref{qeq7}), (\ref{qeq9}) we have the following
theorem.

\begin{theorem}\label{thchar1}
Any quite $o$-minimal Ehrenfeucht theory $T$ has a Rudin--Keisler
preorder, represented by a lattice $Q_k\times L_{s,3}$, and a
decomposition formula of the form
$$
3^k\cdot 6^s=2^k\cdot 3^s+\sum\limits_{t=0}^k\sum\limits_{m=0}^s
2^{s-m}\cdot(2^t\cdot 4^m-1)\cdot C^t_k\cdot C^m_s.
$$
For $s=0$ the decomposition formula has the form {\rm
(\ref{qeq5})}, and for $k=0$~--- {\rm (\ref{qeq7})}.
\end{theorem}

\end{document}